\documentclass[twoside]{article}
\usepackage{amsmath,amsfonts,amssymb,amsthm}
\usepackage{hyperref}
\usepackage[latin2]{inputenc}
\usepackage[a4paper]{geometry}
\usepackage[color,matrix,arrow,all]{xy}
\usepackage{stackrel}
\usepackage{relsize}
\newcommand{\ds}{\displaystyle\sum}
\newcommand{\dip}{\displaystyle\prod}

\newcommand{\complex}{\mathbb{C}}
\newcommand{\nat}{\mathbb{N}}
\newcommand{\integ}{\mathbb{Z}}
\newcommand{\orb}{\mathcal{O}}
\newcommand{\op}{\operatorname}
\newcommand{\Hom}{\operatorname{Hom}}

\newcommand{\Ext}{\operatorname{Ext}}
\newcommand{\GL}{\operatorname{GL}}
\newcommand{\SL}{\operatorname{SL}}
\newcommand{\Rep}{\operatorname{Rep}}
\newtheorem{theorem}{Theorem}[section]

\newtheorem{lemma}[theorem]{Lemma}
\newtheorem{proposition}[theorem]{Proposition}
\newtheorem{corollary}[theorem]{Corollary}

\theoremstyle{definition}
\newtheorem{remark}[theorem]{Remark}

\newtheorem{notation}[theorem]{Notation}

\newtheorem{example}[theorem]{Example}

\def\presuper#1#2%
  {\mathop{}%
   \mathopen{\vphantom{#2}}^{#1}%
   \kern-\scriptspace%
   #2}

\title{The $b$-functions of semi-invariants of quivers}
\author{Andr\'as Cristian L\H{o}rincz}
\date{}

\begin{document}

\maketitle

\begin{abstract}
In this paper we compute $b$-functions (or Bernstein-Sato polynomials) of various semi-invariants of quivers. The main tool is an explicit relation for the $b$-functions between semi-invariants that correspond to each other under reflection functors (or castling transforms). This enables us to compute recursively the $b$-functions for all Dynkin quivers, and extended Dynkin quivers with prehomogeneous dimension vectors. 
\end{abstract}

\section*{Introduction}
\label{sec:intro}

A lot of effort has been made to calculate the $b$-functions of semi-invariants of prehomogeneous vector spaces (for example, see \cite{kimu,skko,ukai}). In \cite{graschu} the $b$-functions of linear free divisors of prehomogeneous spaces are studied, among them some $b$-functions arising from quivers. In \cite{sugi} the $b$-functions for quivers of type $A_n$ are computed. This paper is a natural continuation of the program.

Take $X=(x_{ij})$ an $n\times n$ generic matrix of variables, and $\partial X$ the matrix formed by the partial derivatives $\dfrac{\partial}{\partial x_{ij}}$. Its determinant is a differential operator. The following is Cayley's classical identity:

$$\det \partial X \det X^{s+1} = (s+1)(s+2)\cdots(s+n)\det(X)^s.$$

Here the $b$-function of the determinant is $b(s)=(s+1)(s+2)\cdots (s+n)$. 

Our formulas are generalizations of this identity for determinants of various block matrices (see Remark \ref{block}; see Section \ref{sec:reflect} for examples).

In Section \ref{sec:b-func}, we focus on generalities about $b$-functions, mostly in the equivariant setting. 

In Section \ref{sec:cast} we give a relation of $b$-functions of semi-invariants related under castling transforms (Theorem \ref{thm:castle}).

In Section \ref{sec:quiv}, we give some background material on quivers and their semi-invariants.

In Section \ref{sec:reflect}, we apply the results from Section \ref{sec:cast} in the context of quivers. In particular, the technique allows us to compute using reflection functors all $b$-functions (of one and several variables) for Dynkin quivers, and extended Dynkin quivers with prehomogeneous dimension vectors. We provide several examples of how this can be done (Examples \ref{ex:first},\ref{ex:second},\ref{ex:third}). We also give an example with a quiver of wild type where the method fails to give the $b$-function directly (see Example \ref{ex:not}).

We provide some geometric consequences of these computations in a following paper \cite{lol}.

We mention that in the preprint \cite{lol2}, we provide an additional technique of computation of $b$-functions using a ``slice technique'' that is also applicable to some quivers. However, the method of reflection functors (i.e. castling transforms) is more useful in this context, as it is applicable to a broader class of semi-invariants of quivers.

\vspace{0.05in}

\begin{notation}
As usual, $\nat$ will denote the set of all non-negative integers and $\complex$ the set of complex numbers. 

For $a,b,d\in \nat, a\leq b$, we use the following notation in $\complex[s]$: 
$$[s]^d_{a,b}:=\prod_{i=a+1}^b \prod_{j=0}^{d-1} (ds+i+j).$$
In the case $d=1$, we sometimes write $[s]_{a,b}:=[s]_{a,b}^1$. Also, if $a=0$, we sometimes write $[s]^d_{b}:=[s]^d_{0,b}$. Hence $[s]^{d}_{a,b}[s]^d_{a}=[s]^d_b$.

Now fix an $l$-tuple $\underline{m}=(m_1,\dots,m_l)\in \nat^l$. Then for any $l$-tuple $(d_1,\dots,d_l)$, we use the following notation in $\complex[s_1,\dots,s_l]$:
$$[s]^{d_1,\dots,d_l}_{a,b}=\prod_{i=a+1}^b \prod_{j=0}^{d-1} (d_1s_1+\dots+d_ls_l+i+j),$$
where $d=m_1d_1+\dots + m_ld_l$.
\end{notation}

\section{b-functions}
\label{sec:b-func}

First we define and briefly recall some basic properties of $b$-functions. For details we refer the reader to \cite{gyoja,kashi}. 

Throughout this paper we work over the complex field $\complex$. Let $V$ be an $n$-dimensional vector space. Denote by $D$ the algebra of differential operators on $V$ (i.e. the Weyl algebra in $n$ variables).

Let $f\in \complex[V]$ be a non-zero polynomial. Then there is a differential operator $P(s)\in D[s]:=D\otimes \complex[s]$ and a non-zero polynomial $b(s)\in \complex[s]$ such that 
$$P(s)\cdot f^{s+1}(x) = b(s)\cdot f^s(x).$$
We call such a function $b(s)$ the $b$-\textit{function} (or Bernstein-Sato polynomial) of $f$ if it is monic of minimal degree. By \cite{kashiwara}, all roots of $b(s)$ are negative rational numbers.

For now, let $G$ be a connected reductive algebraic group, acting rationally on $V$. Then we have an action of $G$ on $\complex[V]$ by $(g\cdot f)(v)=f(g^{-1}v)$, for $g\in G$, $f\in \complex[V]$. We call a polynomial $f\in \complex[V]$ a \itshape semi-invariant \normalfont, if there is a character $\sigma\in \Hom(G,\complex^{\times})$ such that $g\cdot f = \sigma(g) f$, that is, $f(gv)=\sigma(g)^{-1}f(v)$, and in this case we say the weight of $f$ is $\sigma$. We note that several authors, including \cite{kashi}, use the inverse convention for weights and the term ``relative invariant'' for semi-invariant.

We form the ring of semi-invariants 

$$\op{SI}(G,V)=\bigoplus_{\sigma} \op{SI}(G,V)_{\sigma},$$
where the sum runs over all characters $\sigma$ and the weight spaces are
$$\op{SI}(G,V)_{\sigma}=\{f\in \complex[V] | f \text{ is a semi-invariant of weight } \sigma\}.$$

The multiplicity of $\sigma$ is $\dim\op{SI}(G,V)_{\sigma}$. We say that $\sigma$ is \itshape multiplicity-free\normalfont, if the multiplicity of $\sigma^k$ is $1$, for any $k\in\nat$.

If $f\in \complex[V]$ is a non-zero semi-invariant of multiplicity-free weight $\sigma$, there exists (up to a constant) a non-zero dual semi-invariant $f^*\in \complex[V^*]$ of weight $\sigma^{-1}$, which we view as a constant coefficient differential operator.

The following result follows by \cite[Corollary 2.5.10]{gyoja}:

\begin{theorem}
Let $f$ be semi-invariant of multiplicity-free weight and take $f^*$ as above. Then we have the following equation
\begin{equation}\label{eq:good}
f^*(\partial x)f(x)^{s+1}=b(s)f(x)^s,
\end{equation}
and $b(s)$ coincides with the $b$-function of $f$ (up to a non-zero constant factor). Moreover, $\deg b(s) = \deg f$.
\end{theorem}

We call $(G,V)$ a \textit{prehomogeneous vector space}, if $V$ has a dense open orbit $\orb$, i.e. $\overline{\orb}=V$. By a standard result (cf. \cite{preh}), $(G,V)$ is prehomogeneous iff all weight multiplicities of the ring of semi-invariants are at most $1$. 

Now we introduce the more general notion of $b$-functions of several variables. For $f_1,\dots f_l$ semi-invariants of weights $\sigma_1,\dots, \sigma_l$, take respective duals $f^*_1,\dots, f^*_l$, and put $\underline{f}=(f_1,\dots, f_l)$ and $\underline{f}^*=(f_1^*,\dots , f_l^*)$. For a multi-variable $\underline{s}=(s_1,\dots ,s_l)$, we define their powers by $\underline{f}^{\underline{s}}=\dip_{i=1}^l f_i^{s_i}$, and $\underline{f}^{*\underline{s}}=\dip_{i=1}^l f_i^{*s_i}$. Also, we say that the collection $\underline{\sigma}=(\sigma_1,\dots,\sigma_l)$ is \textit{multiplicity-free}, if the multiplicity of $\sigma_1^{k_1}\cdots\sigma_l^{k_l}$ is $1$, for any $(k_1,\dots,k_l)\in \nat^l$. This is equivalent to the product $\sigma_1\cdots \sigma_l$ being multiplicity-free.

\begin{lemma}
Using the notation above, if $\underline{\sigma}$ is multiplicity-free, then for any $l$-tuple $\underline{m}=(m_1,\dots,m_l)\in \nat^l$ there is a polynomial $b_{\underline{m}}(\underline{s})$ of $l$ variables such that 
\begin{equation}
\underline{f}^{*\underline{m}}(\partial x)\underline{f}^{\underline{s}+\underline{m}}(x)=b_{\underline{m}}(\underline{s})\underline{f}^{\underline{s}}(x).
\end{equation}
\end{lemma}

If $\underline{\sigma}$ is multiplicity-free, then all weights $\sigma_i$ are multiplicity-free, and one can easily recover the $b$-function $b_{f_i}$ of one variable from $b_{\underline{m}}$. Again, if $(G,V)$ is prehomogeneous, then $\underline{\sigma}$ is multiplicity-free.

\section{b-functions under castling transforms}
\label{sec:cast}

Let $G$ be a connected reductive algebraic group, let $(\pi,V)$ and $(\rho,W)$ be two finite-dimensional rational representations of $G$ and fix $\dim V=n$. Denote by $\Lambda_1$ the standard representation of $\GL$. Take two numbers $r_1,r_2\in\nat$ such that $r_1+r_2=n$. Following \cite[Section 2.3]{kac}, we form two representations:

$$R_1=(G\times \GL_{r_1},(\pi^*\otimes \Lambda_1)\oplus (\rho\otimes 1),V^{r_1}\oplus W),$$
$$R_2=(G\times \GL_{r_2},(\pi\otimes \Lambda_1^*)\oplus (\rho\otimes 1),(V^*)^{r_2}\oplus W).$$

In \cite{saki} such representations $R_1,R_2$ are said to be \itshape castling transforms \normalfont of each other, while in representation theory of quivers the functors relating the representation spaces are called \itshape reflection functors \normalfont. By \cite{kac}, there are canonical isomorphisms of rings of invariants

\begin{equation}\label{eq:isom}
\begin{array}{ll}
\complex[R_1]^{G\times \op{SL}_{r_1}}\cong \complex[R_2]^{G\times\op{SL}_{r_2}}, & \text{ when } r_1,r_2>0, \\[0.2cm] 
\complex[R_1]^{G\times \op{SL}_{r_1}}\cong \complex[R_2]^G\otimes \complex[\op{det}_{r_1}], & \text{ when } r_2=0.
\end{array}
\end{equation}

The papers \cite{kimu,saoc} give relations between the $b$-functions of semi-invariants of prehomogeneous spaces related under castling transform with some extra hypothesis (so-called regularity condition). However, this condition is too restrictive for our purposes, as it is rarely satisfied for quivers. We give an extended result, for the regularity condition turns out to be unnecessary. Additionally, the auxiliary representation $W$ is important, since for quivers we apply castling at only one vertex at a time, while the rest of the quiver will remain unchanged. The proof we give is similar to the sketch of proof of \cite[Theorem 7.51.]{preh}.

Let $f \in \complex[R_1]$ and $f'\in \complex[R_2]$ be two semi-invariants (so $[G,G]\times \op{SL}_{r_i}$-invariants) corresponding under the isomorphisms above. Let $d$ be the absolute value of their $\GL_{r_i}$-weights (they are equal).

\begin{theorem}\label{thm:castle}
Assume $f\in \complex[R_1]$ and $f'\in \complex[R_2]$ are $G\times \GL_{r_i}$-semi-invariants with multiplicity-free weights corresponding under the isomorphisms (\ref{eq:isom}). Then their $b$-functions satisfy
$$b_{f}(s)\cdot [s]^{d}_{r_2} = b_{f'}(s)\cdot [s]^{d}_{r_1}.$$
\end{theorem}
\begin{proof}
The case $r_2=0$ is easy to check directly, since the second isomorphism in (\ref{eq:isom}) gives a separation of variables. So we can assume $r_1,r_2>0$. 
Let $x_{ij}$ $(1\leq i\leq n, 1\leq j \leq r_1)$ and $y_{ij}$ $(1\leq i\leq n,1\leq j\leq r_2)$ be indeterminates in $\complex[V^{r_1}]$ and $\complex[(V^*)^{r_2}]$ respectively. Also, put 
$$\Lambda:=\{ \lambda=(i_1,\dots,i_{r_1}) \; | \; 0<i_1<\dots<i_{r_1}\leq n\},$$
$$\Lambda':=\{ \lambda'=(j_1,\dots,j_{r_2}) \; | \; 0<j_1<\dots<j_{r_2}\leq n\}.$$
For $\lambda=(i_1,\dots,i_{r_1})\in \Lambda$ (resp. for $\lambda'=(j_1,\dots,j_{r_2})\in \Lambda'$), put $|\lambda|=i_1+\dots + i_{r_1}$ (resp. $|\lambda'|=j_1+\dots + j_{r_2}$) and
\[
x_\lambda=\det \begin{pmatrix}
x_{i_1,1} & \dots & x_{i_1,r_1}\\
\vdots & \ddots & \vdots \\
x_{i_{r_1},1} & \dots & x_{i_{r_1},r_1}
\end{pmatrix}
\text {, resp. }
y_{\lambda'}=\det \begin{pmatrix}
y_{j_1,1} & \dots & y_{j_1,r_2}\\
\vdots & \ddots & \vdots \\
y_{j_{r_2},1} & \dots & y_{j_{r_2},r_2}
\end{pmatrix}.
\]

Let $A$ (resp. $A'$) be the subring of $\complex[R_1]$ (resp. $\complex[R_2]$) generated by the polynomials $x_\lambda$, where $\lambda\in\Lambda$ (resp. $y_{\lambda'}$, where $\lambda'\in \Lambda'$).
Let $A_k$ (resp. $A'_k$) denote its homogeneous part of degree $r_1k$ (resp. $r_2k$). Similarly, we define the ring of (constant) differential operators $D$ and $D'$ generated by $\partial x_{\lambda}$ (resp. $\partial y_{\lambda'})$, where
\[
\partial x_\lambda=\det \begin{pmatrix}
\partial x_{i_1,1} & \dots & \partial x_{i_1,r_1}\\
\vdots & \ddots & \vdots \\
\partial x_{i_{r_1},1} & \dots & \partial x_{i_{r_1},r_1}
\end{pmatrix}
\text {, resp. }
\partial y_{\lambda'}=\det \begin{pmatrix}
\partial y_{j_1,1} & \dots & \partial y_{j_1,r_2}\\
\vdots & \ddots & \vdots \\
\partial y_{j_{r_2},1} & \dots & \partial y_{j_{r_2},r_2}
\end{pmatrix}.
\]

Let $D_k$ (resp. $D'_k$) denote the homogeneous part of degree $r_1k$ (resp. $r_2k$).

Now we endow $A$ with the natural action of $\GL_n$. In fact, $A$ can be viewed as the coordinate algebra of the affine Grassmannian $\widetilde{\op{Gr}}(r_1,V)$, that is, the affine cone of the usual Grassmannian variety. Similarly, we equip $A'$ with the dual action of $\GL_n$, viewing it as the coordinate algebra of $\widetilde{\op{Gr}}(r_2,V^*)$. Due to the natural isomorphism $\widetilde{\op{Gr}}(r_1,V)\cong \widetilde{\op{Gr}}(r_2,V^*)$, we have an $\SL_n$-equivariant isomorphism of graded algebras $\tau:A\to A'$. Similarly, we naturally equip $D_k$ (resp. $D'_k$) with the $\GL_n$-structure dual to $A_k$ (resp. $A'_k$) via the pairing $\langle Q(\partial x_\lambda),P(x_\lambda)\rangle=Q(\partial x_\lambda)P(x_\lambda)$.

More explicitly, for $\lambda \in \Lambda$ let $\lambda'\in \Lambda'$ be the complementary set to $\lambda$, namely, $\{\lambda,\lambda'\}=\{1,\dots,n\}$. Then $\tau$ is given by $x_\lambda \mapsto (-1)^{|\lambda'|} y_{\lambda'}$, and $\tau'$ is given by $\partial x_\lambda \mapsto (-1)^{|\lambda'|}\partial y_{\lambda'}$.

For any $k,l\in \nat,k\leq l$,  we have a $\GL_n$-equivariant map $\phi_{k,l}:D_k\otimes A_l \to A_{l-k}$ given by $Q(\partial x_\lambda)\otimes P(x_\lambda)\mapsto Q(\partial x_\lambda)P(x_\lambda)$, and similarly a map $\phi'_{k,l}:D'_k\otimes A'_l\to A'_{l-k}$. So $\phi_{k,l}$ and $\tau^{-1}\circ \phi' _{k,l}\circ (\tau'\otimes \tau)$ are two $\GL_n$-module morphisms $D_k\otimes A_l \to A_{l-k}$:

\[\xymatrix@C+1pc@R+1pc{
D_k\otimes A_l \ar[r]^{\phi_{k,l}} \ar[d]_{\tau'\otimes \tau} & A_{l-k}\\
D'_k\otimes A'_l \ar[r]^{\phi'_{k,l}} & A'_{l-k} \ar[u]_{\tau^{-1}}
}\]

We claim that the diagram commutes up to a constant $c_{k,l}\in \complex$. It is well-known that $A_k$ is an irreducible $\GL_n$-representation corresponding to the Young tableaux of rectangular shape having $r_1$ rows and $l$ columns (see for instance \cite[Proposition 3.1.4]{jerzy}). Using this and the Littlewood-Richardson rule (cf. \cite{jerzy}), one easily sees that in the decomposition of $D_k\otimes A_l$ into irreducible $\GL_n$-modules, $A_{l-k}$ appears with multiplicity $1$. This in turn implies using Schur's lemma that there is a constant $c_{k,l}$ such that $\phi=c_{k,l} \cdot \tau^{-1}\circ \phi' _{k,l}\circ (\tau'\otimes \tau)$. 

To determine $c_{k,l}$, we look on the value of $\phi_{k,l}$ and $\tau^{-1}\circ \phi' _{k,l}\circ (\tau'\otimes \tau)$ on 
$$\partial x_{(1,\dots, r_1)}^k \otimes x^l_{(1,\dots,r_1)}.$$
Using the classical Cayley identity for the determinant, we get
$$\phi_{k,l}(\partial x_{(1,\dots, r_1)}^k \otimes x^l_{(1,\dots,r_1)})= \prod_{i=0}^{k-1} \prod_{j=0}^{r_1-1} (l-i+j) x^{l-k}_{(1,\dots,r_1)},$$
$$\left(\tau^{-1}\circ \phi' _{k,l}\circ (\tau'\otimes \tau)\right)(\partial x_{(1,\dots, r_1)}^k \otimes x^l_{(1,\dots,r_1)})= \prod_{i=0}^{k-1} \prod_{j=0}^{r_2-1} (l-i+j)  x^{l-k}_{(1,\dots,r_1)}.$$
From this, we get $c_{k,l}=\dfrac{\prod_{i=0}^{k-1} \prod_{j=0}^{r_2-1} (l-i+j)}{\prod_{i=0}^{k-1} \prod_{j=0}^{r_1-1} (l-i+j)}$.

\vspace{0.05in}

Now, by the First Fundamental Theorem for $\op{SL}$ (cf. \cite{vinpop}), we know that $f\in A_{d}\otimes \complex[W]$, $f'\in A'_{d}\otimes\complex[W]$, and $(\tau \otimes 1)(f)=f'$, $(\tau'\otimes 1)(f^*)=f'^*$. Moreover, $f^s \in A_{sd} \otimes \complex[W]$, and $(\tau\otimes 1)(f^s)=f'^s\in A'_{sd}\otimes \complex[W]$ for any $s\in\nat$. Separating variables, the computations above yield
$$b_f(s)=c_{d,d(s+1)}\cdot b_{f'}(s),$$
for any $s\in\nat$, hence the conclusion.
\end{proof}

By the same argument, we give the version for the $b$-function of several variables. 

Let $f_i\in \complex[R_1]$ and $f'_i\in \complex[R_2]$, where $i=1,\dots ,l$ be semi-invariants corresponding respectively under the isomorphisms (\ref{eq:isom}), such that the $l$-tuples $\underline{f}=(f_1,\dots,f_l)$ and $\underline{f}'=(f'_1,\dots, f'_l)$ have multiplicity-free weights. Denote by $d_i$ the $\GL$-weight of $f_i$ and $f'_i$.

\begin{theorem}\label{thm:multicastle}
Using the notation above, the $b$-functions of $\underline{f}$ and $\underline{f}'$ satisfy
$$b_{\underline{f}}(s)\cdot [s]^{d_1,\dots,d_l}_{r_2} = b_{\underline{f}'}(s)\cdot [s]^{d_1,\dots,d_l}_{r_1}.$$
\end{theorem}

\section{Semi-invariants of quivers}
\label{sec:quiv}

In this section, we will introduce some basics of quivers and semi-invariants, and we consider slices in this setting. For more on this material, we refer the reader to \cite{elements,brion,harwey}.

A quiver $Q$ is an oriented graph, i.e. a pair $Q=(Q_0,Q_1)$ formed by the set of vertices $Q_0$ and the set of arrows $Q_1$. An arrow $a$ has a head $ha$, and tail $ta$, that are elements in $Q_0$:

\[\xymatrix{
ta \ar[r]^{a} & ha
}\]

A vertex $x\in Q_0$ is called a sink (resp. source) if there is no arrow in $Q$ starting (resp. ending) in $x$. A representation $V$ of $Q$ is a family of finite dimensional vector spaces $\{V(x)\,|\, x\in Q_0\}$ together with linear maps $\{V(a) : V(ta)\to V(ha)\, | \, a\in Q_1\}$. The dimension vector $\underline{d} (V)\in \nat^{Q_0}$ of a representation $V$ is the tuple $\underline{d}(V):=(d_x)_{x\in Q_0}$, with $d_x=\dim V(x)$. A morphism $\phi:V\to W$ of two representations is a collection of linear maps $\phi = \{\phi(x) : V(x) \to W(x)\,| \,x\in Q_0\}$, with the property that for each $a\in Q_1$ we have $\phi(ha)V(a)=W(a)\phi(ta)$. Denote by $\Hom_Q(V,W)$ the vector space of morphisms of representations from $V$ to $W$. For two vectors $\alpha, \beta\in \integ^{Q_0}$, we define the Euler product
$$\langle \alpha, \beta \rangle = \ds_{x\in Q_0} \alpha_x \beta_x - \ds_{a\in Q_1} \alpha_{ta} \beta_{ha}.$$
We define the vector space of representations with dimension vector $\alpha\in \nat^{Q_0}$ by
$$\Rep(Q,\alpha):=\displaystyle\bigoplus_{a\in Q_1} \Hom(\complex^{\alpha_{ta}},\complex^{\alpha_{ha}}).$$
The following group acts on $\Rep(Q,\alpha)$ by conjugation in the obvious way:
$$\GL(\alpha):= \prod_{x\in Q_0} \GL_{\alpha_x}.$$
Under this action, two representations lie in the same orbit iff they are isomorphic representations.

For any two representations $V$ and $W$, we have the following exact sequence:

\begin{equation}\label{eq:ringel}
\begin{array}{rlc}
0 \to \Hom_Q (V,W) \stackrel{i}{\longrightarrow} \displaystyle\bigoplus_{x \in Q_0}& \!\!\!\!\!\Hom(V(x),W(x)) & \\
& \stackrel{d^V_W}{\longrightarrow}  \displaystyle\bigoplus_{a\in Q_1} \Hom(V(ta),W(ha)) \stackrel{p}{\longrightarrow} \Ext_Q(V,W)\to 0 &
\end{array}
\end{equation}

Here, the map $i$ is the inclusion, $d_W^V$ is given by
$$\{\phi(x)\}_{x\in Q_0} \mapsto \{\phi(ha)V(a) - W(a)\phi(ta)\}_{a\in Q_1}$$
and the map $p$ builds an extension of $V$ and $W$ by adding the maps $V(ta)\to W(ha)$ to the direct sum $V\oplus W$.

From the exact sequence (\ref{eq:ringel}) we have that $\langle \underline{d}(V),\underline{d}(W) \rangle = \dim \Hom(V,W) - \dim \Ext(V,W)$.

There is a geometric interpretation of this exact sequence when $V=W$, which gives in particular that the orbit $\mathcal{O}_V$ is dense iff $\Ext_Q(V,V)=0$. In this case we say $V$ is \textit{generic} and $\underline{d}(V)$ is a \textit{prehomogeneous dimension vector}. We call the decomposition $V=\oplus V_i^{m_i}$ of a generic representation into indecomposable representations the canonical decomposition of $\underline{d}(V)$.

Now we turn to semi-invariants of a quiver representation space $\Rep(Q,\beta)$. As in Section \ref{sec:b-func}, form the ring of semi-invariants $\op{SI}(Q,\beta)\subset \complex[\Rep(Q,\beta)]$ by 
$$\op{SI}(Q,\beta)=\bigoplus_\sigma \op{SI}(Q,\beta)_\sigma.$$
Here $\sigma$ runs through all the characters of $\GL(\beta)$. Each character $\sigma$ of $\GL(\beta)$ is a product of determinants, that is, of the form 
$$\prod_{x\in Q_0} \op{det}_x^{\sigma(x)},$$
where $\op{det}_x$ is the determinant function on $\GL_{\beta_x}$. In this way, we will view a character $\sigma$ as a function $\sigma : Q_0 \to \integ$, or equivalently, as an element $\sigma\in \Hom_{\integ}(\integ^{Q_0},\integ)$. With this convention, we view characters as duals to dimension vectors, namely:
$$\sigma(\beta)=\ds_{x\in Q_0} \sigma(x)\beta_x.$$

Throughout we will assume that $Q$ is a quiver \itshape without oriented cycles \normalfont.

We recall the definition of an important class of determinantal semi-invariants, first constructed by Schofield \cite{scho}. Fix two dimension vectors $\alpha,\beta$, such that $\langle \alpha, \beta \rangle=0$. Then for every $V\in \Rep(Q,\alpha)$ and $W\in \Rep(Q,\beta)$ the matrix of the map $d^V_W$ in (\ref{eq:ringel}) will be a square matrix. We define the semi-invariant $c$ of the action of $\GL(\alpha)\times\GL(\beta)$ on $\Rep(Q,\alpha)\times \Rep(Q,\beta)$ by $c(V,W):=\det d^V_W$. For a fixed $V$, restricting $c$ to $\{V\} \times \Rep(Q,\beta)$ defines a semi-invariant $c^V\in \op{SI}(Q,\beta)$. Similarly, for a fixed $W$, restricting $c$ to $\Rep(Q,\alpha)\times \{W\}$, we get a semi-invariant $c_W\in \op{SI}(Q,\alpha)$. The weight of $c^V$ is $\langle \alpha, \cdot \rangle \in \Hom_{\integ}(\integ^{Q_0},\integ)$, and the weight of $c_W$ is $-\langle \cdot, \beta \rangle$. The semi-invariants $c^V$ and $c_W$ are well-defined up to scalar, that is, they don't depend on the isomorphism classes of $V$ and $W$.

\begin{theorem}[{\cite{harwey,schovan}}]\label{thm:span}
For a fixed vector $\beta$, the ring of semi-invariants $\op{SI}(Q,\beta)$ is spanned by the semi-invariants $c^V$, with $\langle \underline{d}(V),\beta \rangle=0$. The analogous result is true for semi-invariants $c_W$.
\end{theorem}

For a more precise statement when $\beta$ is prehomogeneous, see \cite{scho}.

To find multiplicity-free weights for semi-invariants on $\Rep(Q,\beta)$ with $\beta$ not necessarily prehomogeneous, the following reciprocity result is useful:

\begin{lemma}[{\cite[Corollary 1]{harwey}}]
\label{lem:recip}
Let $\alpha$ and $\beta$ be two dimension vectors, with $\langle \alpha,\beta \rangle=0$. Then
$$\dim \op{SI}(Q,\beta)_{\langle \alpha, \cdot \rangle} = \dim \op{SI}(Q,\alpha)_{-\langle \cdot , \beta \rangle}.$$
\end{lemma}

In particular, if $f$ is a non-zero semi-invariant of weight $\langle \alpha , \cdot \rangle$ with $\alpha$ prehomogeneous, then any multiple of $\alpha$ is also prehomogeneous, hence we see that $f=c^V$ has multiplicity-free weight with $V$ generic in $\Rep(Q,\alpha)$.

\begin{remark}
\label{block}
For any $V\in \Rep(Q,\alpha)$, it is easy to write down the semi-invariants $c^V(W)$ explicitly, as determinants of suitable block matrices. Namely, label the rows formed by the blocks with the arrows $a\in Q_1$, and label the columns with the vertices in $Q_0$. Then, for an arrow $a$, we put two block entries in the row of $a$: $I_{\dim V(ta)}\otimes W(a)$ in the column $ta$, and $-V(a)^t \otimes I_{\dim W(ha)}$  in the column $ha$.
\end{remark}

We give the following easy lemma:

\begin{lemma}\label{thm:simp}
Let $Q$ be a quiver without oriented cycles, $\beta$ a dimension vector and $f$ a semi-invariant on $\Rep(Q,\beta)$ of weight $\sigma=\langle \alpha , \cdot \rangle$. Then we can view $f$ as a semi-invariant on a new quiver with new weight according to the following simplification rules:
\begin{itemize}
\item[(a)] If $\alpha_1=0$, then we have (we put the values of $\alpha$ on top of $\beta$):
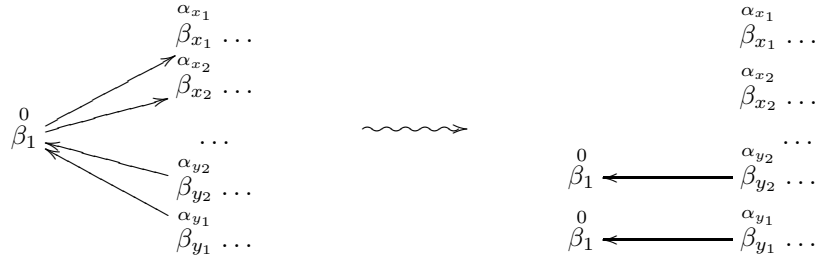
 \[\vcenter{\vbox{\xymatrix@R-2.3pc@C+1.8pc{
 & \stackrel{\alpha_{x_1}}{\beta_{x_1}} \dots \\
 & \stackrel{\alpha_{x_2}}{\beta_{x_2}} \dots \\
\stackrel{0}{\beta_{1}} \ar[uur] \ar[ur] & \dots \\
 & \stackrel{\alpha_{y_2}}{\beta_{y_2}} \ar[ul] \dots \\
 & \stackrel{\alpha_{y_1}}{\beta_{y_1}} \ar[uul] \dots
}}}\hspace{0.4in}{\xymatrix@C+1.5pc{ \ar@{~>}[r] & }}\hspace{0.4in}
 \vcenter{\vbox{\xymatrix@R-2pc@C+2pc{
 & \stackrel{\alpha_{x_1}}{\beta_{x_1}} \dots \\
 & \stackrel{\alpha_{x_2}}{\beta_{x_2}} \dots \\
 & \hspace{0.2in} \dots \\
\stackrel{0}{\beta_{1}} & \stackrel{\alpha_{y_2}}{\beta_{y_2}} \ar[l] \dots \\
\stackrel{0}{\beta_{1}} & \stackrel{\alpha_{y_1}}{\beta_{y_1}} \ar[l] \dots
}}} \]
\item[(b)] Write $\sigma=-\langle \cdot, \alpha^* \rangle$. If $\alpha^*_1=0$, then the same simplification rule holds as in part (a) by replacing $\alpha$ with $\alpha^*$, with the arrows reversed.
\end{itemize}
\end{lemma}

\begin{proof}
\begin{itemize}
\item[(a)] We can assume $f=c^V$, with $\underline{d}(V)=\alpha$. Then we see explicitly that $f$ doesn't depend on the arrows from $1$ to $x_i$, hence we can drop them. Finally, we can split vertex $1$ so that the arrows from $y_i$ have different heads, not changing $f$.
\item[(b)] We can assume $f=c_W$, with $\underline{d}(W)=\gamma$. Then the we see the simplifications explicitly as in part (a).
\end{itemize}
\end{proof}

\section{b-functions of semi-invariants of quivers}
\label{sec:reflect}

First, we introduce some terminology for reflection functors, that is, castling transforms in the quiver setting (for details, see \cite[VII.5.]{elements}). Throughout, let $Q$ be a quiver without oriented cycles. Given any vertex $i\in Q_0$, we form a new quiver $c_i Q$ by reversing all arrows that start or end in $i$. An ordering of $i_1,\dots,i_n$ of the vertices of $Q$ is called admissible if for each $k$ the vertex $i_k$ is a sink for $c_{i_{k-1}}\dots c_{i_1} Q$. In such case, it is easy to see that $c_{i_n}\dots c_{i_1}Q=Q$. Since $Q$ has no oriented cycles, $Q$ has admissible orderings, and we fix one. For $x\in Q_0$, we take the following the linear map of dimension vectors that we denote by the same letter
$$c_x:\integ^n\to \integ^n$$
$$c_x(\beta)_y=\begin{cases}
\beta_y &\text{ if } x\neq y,\\
-\beta_x+\ds_{\text{edges } x \text{---} z}\beta_z& \text{ if } x=y.
\end{cases}$$

Also, let $c=c_{i_n}\cdots c_{i_{1}}$ be the Coxeter transformation. It is independent on the choice of the admissible ordering. As a matrix, we have that $c=-E^{-1} E^t$, where $E$ denotes the Euler matrix of $Q$ corresponding to the Euler product.

The fact that the reflection functors are particular cases of castling transforms has been noted in \cite{kac}. Namely, let $x$ be a sink with $c_x(\beta)_x\geq 0$. Then using the notation from Section \ref{sec:cast}, we put
$$r_1 = \beta_x\,,\, r_2=c_x(\beta)_x \,,\, G=\prod_{y\in Q_0, y\neq x}\GL_{\beta_y} \,,\, V^{r_1}= \bigoplus_{a\in Q_1, ta=x} \Hom(\complex^{\beta_{ha}},\complex^{\beta_x}),$$ 
$$W=\bigoplus_{a\in Q_1, ta\neq x} \Hom(\complex^{\beta_{ha}},\complex^{\beta_{ta}})\,,\, (V^*)^{r_2}=\bigoplus_{a\in c_x(Q)_1, ha=x} \Hom(\complex^{c_x(\beta)_x},\complex^{\beta_{ta}}).$$

The isomorphisms from (\ref{eq:isom}) translate into the quiver setting as:

\begin{equation}\label{eq:isom2}
\begin{array}{ll}
c_x:\op{SI}(Q,\beta)\cong \op{SI}(c_x Q, c_x(\beta)), & \text{ when } c_x(\beta)_x>0, \\[0.2cm] 
c_x:\op{SI}(Q,\beta)\cong \op{SI}(c_x Q, c_x(\beta))\otimes \complex[\op{det}_{\beta_x}], & \text{ when } c_x(\beta)_x=0.\\[0.2cm]
\end{array}
\end{equation}

Note that $\op{SI}(Q,\beta)\cong \op{SI}(Q, \beta-\beta_x \epsilon_x)$, when $c_x(\beta)_x<0$, where $\epsilon_x$ is the dimension vector given by $(\epsilon_x)_x=1$ and $(\epsilon_x)_y=0$ for $y\in Q_0\backslash\{x\}$. This is because in this case a semi-invariant doesn't depend on the arrows ending in $x$, hence we can simply ``drop'' the vertex.  

All these isomorphisms respect weight spaces: when $c_x(\beta)_x\geq 0$, we have 
$$\op{SI}(Q,\beta)_{\langle \alpha,\cdot \rangle}\cong \op{SI}(c_x Q, c_x(\beta))_{\langle c_x(\alpha),\cdot \rangle}.$$

We have the reflection functors on the representation level $C_x:\op{rep}(Q)\to \op{rep}(c_x Q)$ such that $C_x(S_x)=0$, for the simple module $S_x$ at $x$, and for all other indecomposables $X$, $C_x(X)$ is a non-zero indecomposable representation with dimension vector $c_x(\underline{d}(X))$ (see \cite{bgp}). 

Now denote by $C$ the Coxeter functor defined by $C=C_{i_n}\cdots C_{i_{1}}:\op{rep}(Q)\to \op{rep}(Q)$. Then $C(P_y)=0$, for the projective module $P_y$ corresponding to any vertex $y\in Q_0$, and for all other indecomposables $X$, $C(X)$ is a non-zero indecomposable representation with dimension vector $c(\underline{d}(X))$. We say an indecomposable $X$ is preprojective, if $C^k(X)=0$, for some $k\in \nat$. Preprojective representations are generic and we also call their dimension vectors preprojective. Dually, we can define these notions with $x$ a source, we get the Coxeter transformation $c^{-1}$ and preinjective representations. An indecomposable representation is called regular is it is neither preprojective nor preinjective. 

We recall (cf. \cite{elements}) the notions of a Dynkin quiver (of type $A_n, D_n, E_6, E_7, E_8$) and of an extended Dynkin quiver (of type $\tilde{A}_n, \tilde{D}_n, \tilde{E}_6,\tilde{E}_7,\tilde{E}_8$). Since we are dealing with quivers without oriented cycles, we will exclude the cycle quiver of type $\tilde{A}_n$.
By Gabriel's theorem (see \cite{elements}), the quiver is of Dynkin type if and only if all the indecomposables are preprojective (and preinjective).

\begin{theorem}\label{thm:bref} Let $Q$ be a quiver without oriented cycles and $f_i\in \op{SI}(Q,\beta)_{\langle \alpha^i, \cdot \rangle}$ be semi-invariants, where $i=1,\dots,k$. Assume $\underline{f}=(f_1,\dots,f_k)$ has multiplicity-free weight and the coordinates of $c(\beta)$ are non-negative. Then the $b$-function satisfies the formula
$$\mathlarger{b_{\underline{f}}(s)=b_{c(\underline{f})}(s)\mathlarger{\prod_{x\in Q_0}} \dfrac{[s]^{c(\alpha^1)_x,\dots,c(\alpha^k)_x}_{\beta_x}}{[s]^{c(\alpha^1)_x,\dots,c(\alpha^k)_x}_{c(\beta)_x}}}.$$
\end{theorem}

\begin{proof}
Fix a sink $x$. First, note that the case $c_x(\beta)_x< 0$ implies just that none of the semi-invariants depend on $x$ (and all of them have weight $0$ at $x$), so we can drop the vertex. 

Since $x$ is a source in $c_x Q$, the absolute value of the $\GL_{\beta_x}$-weight of $f$ and $c_x(f)$ is $c_x(\alpha^i)_x$. Now applying Theorem \ref{thm:multicastle}, we get that 
$$b_{\underline{f}}(s)=b_{c_x(\underline{f})}(s)\cdot \dfrac{[s]^{c_x(\alpha^1)_x,\dots,c_x(\alpha^k)_x}_{\beta_x}}{[s]^{c_x(\alpha^1)_x,\dots,c_x(\alpha^k)_x}_{c_x(\beta)_x}}.$$
Applying this to an admissible sequence we get the desired formula.
\end{proof}

We call a semi-invariant $f\in \op{SI}(Q,\beta)$ \itshape reducible by reflections \normalfont if, after applying reflection functors finitely many times, we can reduce it to a constant function via the isomorphisms in (\ref{eq:isom2}). In this case, we can compute the $b$-function of $f$ by Theorem \ref{thm:bref}. Similarly, if $f_1,\dots,f_k\in \op{SI}(Q,\beta)$ are all reducible by reflections, then we can compute their $b$-function of several variables.

\begin{proposition}\label{thm:pre}
Let $Q$ be a quiver without oriented cycles, $f\in \op{SI}(Q,\beta)$ a semi-invariant of weight $\langle \alpha,\cdot \rangle$. Assume one of the following cases holds:
\begin{itemize}
\item[(a)] The dimension vector $\alpha$ is preprojective or preinjective, or
\item[(b)] The dimension vector $\beta$ is prehomogeneous and any indecomposable in the canonical decomposition of $\beta$ is preprojective or preinjective.
\end{itemize}
Then $f$ is reducible by reflections.
\end{proposition}

\begin{proof}
\begin{itemize}
\item[(a)] We show this for $\alpha$ preprojective. By Lemma \ref{lem:recip}, $f$ has multiplicity-free weight. Applying the Coxeter transformation sufficiently many times, we arrive at a semi-invariant whose weight corresponds to a projective dimension vector, which implies that it is constant. The only thing one has to deal with is that after applying a reflection transformation $c_x$, one might end up with a vector $\beta'$ with $\beta'_x<0$. But in this case the function doesn't depend on $x$, so after replacing $\beta'_x$ with $0$, we can carry on with the procedure.
\item[(b)] If the canonical decomposition of $\beta$ doesn't have the simple $S_x$ as summand, for $x$ a sink, then applying $c_x$ to each indecomposable in the canonical decomposition gives us the canonical decomposition of $c_x(\beta)$.  We apply the reflection functor in the order given by the admissible sequence repeatedly until we reach a simple, that is, by applying $c_{x_{i-1}}$ we reach the simple $S_{x_{i}}$ as a summand of the canonical decomposition. Write $f=c_W$. Since $f$ doesn't vanish on the generic element, we have that $W_x=\Hom(S_x,W)=0$. By Lemma \ref{thm:simp}, this implies that $f$ doesn't depend on the arrows of the sink $x_i$. So we can drop all $S_{x_{i}}$ from the canonical decomposition and then continue by applying $c_{x_i}$. After we get rid of all preprojectives this way, we start working dually (with sources) to get rid of all preinjectives.
\end{itemize}
\end{proof}

The following corollary follows immediately by either case in Proposition \ref{thm:pre}:

\begin{corollary} All semi-invariants of Dynkin quivers are reducible by reflections.
\end{corollary}

\begin{remark}
The fact that the $b$-function is truly a polynomial in the end is not obvious from the formula given in Theorem \ref{thm:bref}. Also, in case of Dynkin quivers, we can apply the procedure in both directions, either with sinks or with sources.  The $b$-function holds interesting combinatorial information.
\end{remark}

\begin{proposition}
Let $Q$ be an extended Dynkin quiver, and $\beta$ a prehomogeneous dimension vector. Then all semi-invariants in $\op{SI}(Q,\beta)$ are reducible by reflections.
\end{proposition}

\begin{proof}
Using the procedure described in Proposition \ref{thm:pre}, we can reduce to the case when all the indecomposables in the canonical decomposition are regular. By \cite[Lemma 5.1]{zwara}, the left orthogonal category of the generic representation contains a preprojective representation. Hence applying the Coxeter functor sufficiently many times, we arrive at a dimension vector $\beta'$ such that $0=\langle \underline{d}(P_x), \beta' \rangle=\beta'_x $, where $P_x$ is the projective corresponding to a vertex $x$. Hence we can drop vertex $x$, reducing to the Dynkin case.
\end{proof}

Now we consider some examples:

\begin{example}\label{ex:first} We take the following quiver of type $E_6$, $f$ the semi-invariant of weight $\langle \alpha, \cdot \rangle$, we write the values $\alpha$ above the values $\beta$, and $\beta_1+\beta_2+\beta_3+\beta_4+2\beta_5=3\beta_6$, together with the necessary inequalities that assure $f$ is non-zero (see below). The Coxeter transformation is given by
$$c=\begin{pmatrix}
0 & 0 & 0 & 1 & 1 & -1\\
1 & 0 & 0 & 1 & 1 & -1\\
0 & 1 & 0 & 0 & 1 & -1\\
0 & 1 & 1 & 0 & 1 & -1\\
0 & 1 & 0 & 1 & 0 & -1\\
0 & 1 & 0 & 1 & 1 & -1
\end{pmatrix}$$
\[
\vcenter{\vbox{\xymatrix@C-0.7pc@R-1pc{
 & & \stackrel{2}{\beta_5}\ar[d] & & \\
\stackrel{1}{\beta_1} \ar[r]& \stackrel{2}{\beta_2}\ar[r] & \stackrel{3}{\beta_6} & \ar[l] \stackrel{2}{\beta_4}& \ar[l] \stackrel{1}{\beta_3}
}}}\stackrel{\text{Apply } c}{\xymatrix@C+1pc{ \ar@{~>}[r] & }}
\]
\[\vcenter{\vbox{\xymatrix@C-0.5pc@R-1pc{
 & & \stackrel[\beta_2+\beta_4-\beta_6]{}{1}\ar[d] & & \\
\stackrel[\beta_4+\beta_5-\beta_6]{}{1} \ar[r]& \stackrel[\beta_1+\beta_4+\beta_5-\beta_6]{}{2}\ar[r] & \stackrel[\beta_2+\beta_4+\beta_5-\beta_6]{}{3} & \ar[l] \stackrel[\beta_2+\beta_3+\beta_5-\beta_6]{}{2}& \ar[l] \stackrel[\beta_2+\beta_5-\beta_6]{}{1}
}}}
\]
\vspace{0.2in}
\[\vcenter{\vbox{\xymatrix@C-0.5pc@R-1pc{
 & & \stackrel[\beta_1+\beta_3+\beta_5-\beta_6]{}{1}\ar[d] & & \\
\stackrel[\beta_2+\beta_3-\beta_6]{}{0} \ar[r]& \stackrel[\beta_6-\beta_5-\beta_1]{}{1}\ar[r] & \stackrel[\beta_6-\beta_5]{}{2} & \ar[l] \stackrel[\beta_6-\beta_5-\beta_3]{}{1}& \ar[l] \stackrel[\beta_1+\beta_4-\beta_6]{}{0}
}}}
\]

Here, using the Lemma \ref{thm:simp} (a), we can simplify by dropping vertex $1$ and $3$, and as the last step we apply only $c_6$:

\[\vcenter{\vbox{\xymatrix@R-0.5pc{
 & \stackrel[\beta_1+\beta_3+\beta_5-\beta_6]{}{1} & \\
 \stackrel[\beta_6-\beta_5-\beta_1]{}{1} & \ar[l] \stackrel[0]{}{1} \ar[r]  \ar[u]& \stackrel[\beta_6-\beta_5-\beta_3]{}{1}
}}}
\]

Here the semi-invariant becomes constant, so we stop. Hence the $b$-function is

$$b(s)=\Bigg(\dfrac{[s]_{\beta_1}\cdot[s]^2_{\beta_2}\cdot[s]_{\beta_3}\cdot[s]^2_{\beta_4}\cdot[s]_{\beta_5}\cdot[s]^3_{\beta_6}}{[s]_{\beta_4+\beta_5-\beta_6}[s]^2_{\beta_1+\beta_4+\beta_5-\beta_6}[s]_{\beta_2+\beta_5-\beta_6}[s]^2_{\beta_2+\beta_3+\beta_5-\beta_6}[s]_{\beta_2+\beta_4-\beta_6}[s]^3_{\beta_2+\beta_4+\beta_5-\beta_6}}  \Bigg)\cdot$$ $$\cdot\Bigg(\dfrac{[s]_{\beta_1+\beta_4+\beta_5-\beta_6}[s]_{\beta_2+\beta_3+\beta_5-\beta_6}[s]_{\beta_2+\beta_4-\beta_6}[s]^2_{\beta_2+\beta_4+\beta_5-\beta_6}}{[s]_{\beta_6-\beta_5-\beta_1}[s]_{\beta_6-\beta_5-\beta_3}[s]_{\beta_1+\beta_3+\beta_5-\beta_6}[s]^2_{\beta_6-\beta_5}} \Bigg) \cdot [s]_{\beta_6-\beta_5}$$
Note that using the inequalities between $\beta_i$, we can reduce this expression to a polynomial. We can read the inequalities looking at the dimension vectors in the Coxeter transformations:
$$\beta_6\leq \beta_4+\beta_5, \beta_2+\beta_5, \beta_2+\beta_4, \beta_1+\beta_3+\beta_5$$
$$\beta_6\geq \beta_1+\beta_5,\beta_3+\beta_5$$
Using these, one way to write the $b$-function as a polynomial is:
$$b(s)=
[s]_{\beta_4+\beta_5-\beta_6,\beta_1+\beta_4+\beta_5-\beta_6}
[s]_{\beta_2+\beta_5-\beta_6,\beta_2+\beta_3+\beta_5-\beta_6}
[s]^2_{\beta_2+\beta_3+\beta_5-\beta_6,\beta_2}
[s]^2_{\beta_1+\beta_4+\beta_5-\beta_6,\beta_4}\cdot$$
$$\cdot[s]^3_{\beta_2+\beta_4+\beta_5-\beta_6,\beta_6}
[s]_{\beta_1+\beta_3+\beta_5-\beta_6,\beta_1}
[s]_{\beta_6-\beta_5-\beta_1,\beta_3}
[s]_{\beta_6-\beta_5-\beta_3,\beta_6-\beta_5}
[s]^2_{\beta_6-\beta_5,\beta_2+\beta_4+\beta_5-\beta_6}
[s]_{\beta_5}
$$
\end{example}

\vspace{0.1in}

In the next example, we compute the $b$-function of several variables of $4$ semi-invariants:

\begin{example}\label{ex:second} Take the following $D_5$ quiver

\[\vcenter{\vbox{\xymatrix{
 & 5 & & \\
1 \ar[r] & 4\ar[u] & 3\ar[l] & \ar[l] 2
}}}\]

with dimension vector $\beta=(n,n,2n,2n,n)$, where $n\in \nat$. There are $4$ fundamental invariants $f_i$, $i=1,\dots,4$, with weights $\langle \alpha^i,\cdot \rangle$, where $\alpha^1=(0,0,1,0,0),\alpha^2=(1,0,0,1,0),\alpha^3=(0,1,1,1,0),\alpha^4=(1,1,1,1,1).$ Explicitly, if we label the generic matrices as
\[\vcenter{\vbox{\xymatrix{
 & n & & \\
n \ar[r]^A & 2n\ar[u]^D & 2n\ar[l]_C & \ar[l]_B n
}}}\]
then the semi-invariants are $f_1=\det C$, $f_2=\det(D A)$, $f_3=\det(D C B)$, $f_4=\det(A,C B)$. The Coxeter transformation is given by

$$c=\begin{pmatrix}
 0 & 0 & 1 & 0 & -1\\
 1 & 0 & 0 & 0 & -1\\
 1 & 1 & 0 & 0 & -1\\
 1 & 0 & 1 & 0 & -1\\
 0 & 0 & 0 & 1 & -1
\end{pmatrix}$$

Again, we write the (ordered) values of $\alpha$ on top of $\beta$, and apply $c$:
\[\vcenter{\vbox{\xymatrix{
 & \stackrel{0,0,0,1}{n} & & \\
\stackrel{0,1,0,1}{n} \ar[r] & \stackrel{0,1,1,1}{2n}\ar[u] & \stackrel{1,0,1,1}{2n}\ar[l] & \ar[l] \stackrel{0,0,1,1}{n}
}}}\stackrel{\text{Apply } c}{\xymatrix@C+1pc{ \ar@{~>}[r] & }}\]
\[\vcenter{\vbox{\xymatrix{
 & \stackrel{0,1,1,0}{n} & & \\
\stackrel{1,0,1,0}{n} \ar[r] & \stackrel{1,1,1,1}{2n}\ar[u] & \stackrel{0,1,1,1}{n}\ar[l] & \ar[l] \stackrel{0,1,0,0}{0}
}}}\]
Here we can drop vertex $2$, and the second semi-invariant, since it became constant. So we are left with
\[\vcenter{\vbox{\xymatrix{
 & \stackrel{0,1,0}{n} & \\
\stackrel{1,1,0}{n} \ar[r] & \stackrel{1,1,1}{2n}\ar[u] & \stackrel{0,1,1}{n}\ar[l]
}}}\stackrel{}{\xymatrix@C+1pc{ \ar@{~>}[r] & }}\vcenter{\vbox{\xymatrix{
 & \stackrel{1,0,1}{n} & \\
\stackrel{0,0,1}{0} \ar[r] & \stackrel{1,1,1}{n}\ar[u] & \stackrel{1,0,0}{0}\ar[l]
}}}\]
Here the first and last semi-invariants became constant, and the middle is the determinant of size $n$. So the $b$-function is:
$$b_{\underline{f}}(s)=[s]_n^{0,1,0,0}\cdot[s]^{0,1,1,1}_{n,2n}\cdot[s]_n^{0,0,0,1}\cdot[s]_n^{1,0,0,0}\cdot[s]^{1,0,1,1}_{n,2n}\cdot [s]_n^{0,0,1,0}.$$
\end{example}

\vspace{0.1in}

\begin{example}\label{ex:third} Let $Q$ be the Kronecker quiver
\[\xymatrix{ 1 \ar@<0.5ex>[r]\ar@<-0.5ex>[r] & 2}\]
with the prehomogeneous dimension vector $\beta=(n\cdot k, (n+1)k)$. We have the following semi-invariant $f_n$ with weight $\langle \alpha, \cdot \rangle$ corresponding to $\alpha=(n+1,n+2)$:
$$f_n=\det\begin{pmatrix}X&Y&0&\dots&0&0\\0& X & Y &\dots & 0 & 0 \\ &\vdots &&\ddots &&\\ 0 & 0 & 0 &\dots & Y & 0\\ 0 & 0 & 0 &\dots&X&Y\end{pmatrix}.$$
Here there are $n+1$ block columns and $n$ block rows. Applying $c_2$ we get the same quiver (after renumbering) with dimension vector $c_2(\beta)=((n-1)k,nk)$ and $c_2(\alpha)=(n,n+1)$. Hence
$$b_{f_n}(s)=b_{f_{n-1}}\cdot [s]^n_{(n-1)k,(n+1)k}=\prod_{i=1}^n [s]^i_{(i-1)k,(i+1)k}.$$
\end{example}

In the following example (compare with example in \cite{zwara}), we see that the method by reflection functors does not give the full $b$-function of several variables. We outline a potential fix for this:

\begin{example}\label{ex:not} Let $(Q,\beta)$ be the following quiver
\[\xymatrix{
8n \ar[drr] & 3n \ar[dr] & 3n \ar[d] & 3n\ar[dl] & 3n\ar[dll]\\
& & 9n & &
}\]
where $n\in \nat$. Then $\beta$ is a prehomogeneous dimension vector, and is a multiple of the regular indecomposable $\arraycolsep=1.5pt\begin{array}{ccccc} 8 & 3 & 3 & 3 & 3\\ &&9&& \end{array}$ . We have $5$ fundamental semi-invariants $f_i=c^{V_i}\in \op{SI}(Q,\beta),i=1,\dots,5$, where the $\underline{d}(V_i)=\alpha_i$ correspond to the following simples in the left orthogonal category:
\[\arraycolsep=1.5pt
\begin{array}{ccccc}
0 & 0 & 1 & 1 & 1\\
&&2&&
\end{array}\,,\,
\begin{array}{ccccc}
0 & 1 & 0 & 1 & 1\\
&&2&&
\end{array}\,,\,
\begin{array}{ccccc}
0 & 1 & 1 & 0 & 1\\
&&2&&
\end{array}\,,\,
\begin{array}{ccccc}
0 & 1 & 1 & 1 & 0\\
&&2&&
\end{array}\,,\,
\begin{array}{ccccc}
3 & 1 & 1 & 1 & 1\\
&&3&&
\end{array}.
\]
\end{example}

These are also regular representations. So we see that the method of reflections won't work directly if we want to compute their $b$-function of several variables. However, if we consider only the first $4$ semi-invariants $f_1,f_2,f_3,f_4$, then we can drop the first source vertex (since the semi-invariants have weight $0$), and the resulting quiver being extended Dynkin, the $b$-function $b_{f_1,f_2,f_3,f_4}$ (of $4$ variables) can be computed. The $b$-function (of one variable) of $f_5$ can be computed by applying the inverse Coxeter transformation $c^{-1}$ first, when we get $$c^{-1}(\alpha_5)=\arraycolsep=1.5pt\begin{array}{ccccc}
0 & 2 & 2 & 2 & 2\\
&&5&&
\end{array}.$$
Then we can again drop the first source vertex, and we are left with an extended Dynkin quiver.

Hence, we can compute $b_{f_1,f_2,f_3,f_4}$ and $b_{f_5}$. Now, in order to obtain information about the $b$-function of $5$ variables, one should compute the $a$-function first, and then employ the structure theorem on $b$-functions (\cite[Theorem 2]{sata} or \cite[Theorem 1.3.5]{ukai}, see also \cite{sugi}).

\vspace{0.05in}

\section*{Acknowledgement}
The author is indebted to his advisor Jerzy Weyman for many valuable discussions.

\vspace{0.05in}

\bibliographystyle{amsplain}
\bibliography{biblo}

\end{document}